\documentclass[11pt,leqno]{article}
\usepackage[a4paper, total={6in, 8in}]{geometry}
\usepackage[T1]{fontenc}
\usepackage{amssymb}
\usepackage{amsmath}
\usepackage[title, toc]{appendix}
\usepackage{amsthm}  
\usepackage{tikz-cd}  
\usepackage{tikz,stackengine}  
\usepackage{enumerate} 
\usepackage{enumitem}  
\usepackage{pdfpages}  
\usepackage{hyperref}
\usepackage{wasysym}  
\usepackage{mathtools} 
\usepackage{stmaryrd} 
\usepackage{graphicx}
\usepackage{xcolor}
\usepackage{mathrsfs}
\usepackage{titlesec}
\usepackage{ytableau}
\usepackage{bbold}
\usepackage[matrix,arrow,curve]{xy}

\hypersetup{
	colorlinks,
	linkcolor={blue!50!black},
	citecolor={blue!50!black},
	urlcolor={blue!80!black}
}
\makeatletter
\def\th@remark{%
	\thm@headfont{\bfseries}%
	\normalfont % body font
	\thm@preskip\topsep \divide\thm@preskip\tw@
	\thm@postskip\thm@preskip
}
\makeatother
\DeclareSymbolFont{extraitalic} {U}{zavm}{m}{it}
\DeclareMathSymbol{\stigma}{\mathord}{extraitalic}{168}
\newtheorem*{proposition*}{Proposition}
\newtheorem*{theorem*}{Theorem}
\newtheorem{theorem}{Theorem}[subsection]
\newtheorem{proposition}[theorem]{Proposition}

\newtheorem{lemma}[theorem]{Lemma}
\newtheorem{corollary}[theorem]{Corollary}%%[section]
\theoremstyle{remark}
\newtheorem{remark}[theorem]{Remark}%%[section]

\makeatletter
\newcommand*\leftdash{\rotatebox[origin=c]{-45}{$\dabar@\dabar@\dabar@$}}
\newcommand*\rightdash{\rotatebox[origin=c]{45}{$\dabar@\dabar@\dabar@$}}
\newcommand*\mondash[1]{\rotatebox[origin=c]{45}{$\dabar@\dabar@\dabar@$}{\kern-1.5ex\raisebox{-.5ex}{}_{#1}}\kern.5ex}
\makeatother
\newcommand{\Qlbar}{\overline{\mathbb{Q}}_{\ell}}

\newcommand{\fl}{\mathcal{B}}

\newcommand{\tilX}{\widetilde{X}}
\newcommand{\cX}{\mathcal{X}}
\newcommand{\F}{\mathcal{F}}

\newcommand{\cO}{\mathcal{O}}
\newcommand{\hc}{\mathfrak{hc}}
\newcommand{\cc}{\underline{\mathbb{k}}}

\newcommand{\aw}{\mathtt{a}_{w_1, w_2}}
\newcommand{\CS}{\mathfrak{C}}

\newcommand{\Av}{\mathrm{Av}}
\newcommand{\For}{\mathrm{For}}

\newcommand{\Id}{\operatorname{Id}}

\newcommand{\Lim}{\mathbb{L}}
\newcommand{\Limg}{\operatorname{\mathbb{L}_{\gamma}}}
\newcommand{\Hf}{\mathbb{H}}
\newcommand{\Hfb}{\mathbb{H}_{\bullet}}

\newcommand{\Ad}{\operatorname{Ad}}

\titleformat{\section}
{\normalfont\bfseries}{\thesection}{1em}{}

\titleformat{\subsection}[runin]
{\normalfont\bfseries}{\thesubsection}{1em}{}

\titleformat{\subsubsection}[runin]
{\normalfont\bfseries}{\thesubsubsection}{1em}{}

\newcommand{\Addresses}{{% additional braces for segregating \footnotesize
		\bigskip
		\footnotesize
		R.~Gonin, \textsc{School of Mathematics,
			University of Cardiff, Cardiff, United Kingdom}\par\nopagebreak
		\textit{E-mail address}: \texttt{roman.ro.gonin@gmail.com}
		
		\medskip
		
		A.~Ionov, \textsc{Department of Mathematics, University of
Texas at Austin, Austin, Texas, USA}\par\nopagebreak
		\textit{E-mail address}: \texttt{andrei.ionov@austin.utexas.edu}
		
		\medskip
		
		K.~Tolmachov, \textsc{Department of Mathematics, University of Hamburg,
			Hamburg, Germany}\par\nopagebreak
		\textit{E-mail address}: \texttt{tolmak@khtos.com}
		
}}

\begin{document} \title{On a t-exactness property of the Harish-Chandra transform}

\author{Roman Gonin, Andrei Ionov, Kostiantyn Tolmachov}\date{}
\maketitle \abstract{Using hyperbolic localization, we identify the nearby cycles along the Vinberg de\-ge\-ne\-ra\-tion with the composition of Radon and Harish-Chandra functors, both considered for the category of character sheaves. This provides a new, simple proof of the exactness of this composition, extending previously known results to arbitrary monodromy and more general sheaf-theoretic settings.}

\section{Introduction}
\subsection{Notation.}
For a stack $\cX$ defined over $\mathbb{C}$, let $D^b(\cX)$ be the bounded algebraically constructible derived category of sheaves on $\cX$ with analytic topology, with coefficients in a field $\mathbb{k}$. All the stacks we encounter will be quotient stacks of the form $X/H$, where $X$ is a scheme and $H$ is an algebraic group over $\mathbb{C}$. The sheaf theory for such stacks is described in \cite{bernsteinEquivariantSheavesFunctors2006}. We write $\operatorname{pt} = \operatorname{Spec}(\mathbb{C})$ and denote by $\cc_{\cX}$ the constant sheaf on the stack~$\cX$.  

For an algebraic group $H$ acting on a scheme $X$ and its subgroup $H'$ let 
\[\For^H_{H'}\colon D^b(H\backslash X) \to D^b(H'\backslash X)\] 
denote the forgetful functor, namely, the $*$-pullback along the projection $H' \backslash X \to H \backslash X$. Let 
\[\text{$\Av^H_{H',!}$, $\Av_{H',*}^H$} \colon D^b(H'\backslash X) \to D^b(H\backslash X)\]
stand for the $!$ and $*$-pushforward from $H'\backslash X$ to $H\backslash X$, respectively. %If $H'$ is clear from the context, we will omit it from the notation, writing $\For^H, \Av_!^H, \Av_*^H$, respectively. 

We will adopt the same notation for the right action of the group, writing $\prescript{R}{}\For, \prescript{R}{}\Av$ for the right forgetful and averaging functors, when emphasis is needed.

\subsection{Background and the main result.}\label{subsec:intro_results}
Let $G$ be a reductive group, $B$ a fixed Borel subgroup, $U \subset B$ its unipotent radical, and $T \subset B$ a maximal torus. The Harish-Chandra transform $\hc_!$ for constructible sheaves was introduced by Lusztig \cite{lusztigCharacterSheaves1985} in his study of character sheaves. It is given by a pull-push functor along the horocycle correspondence
\[G/_{\Ad}G \xleftarrow{} G/_{\Ad}B \xrightarrow{} (U\backslash G/U)/_{\Ad}T ,\]
\[\hc_!\colon D^b(G/_{\Ad}G) \to D^b \left( (U\backslash G/U)/_{\Ad}T \right).\]
%\[\hc_! := \Av^U_!\circ \For^G.\]
Here the subscript $!$ in $\hc_!$ indicates that the functor is defined using $!$-averaging.

It was observed by Ginzburg \cite{ginzburgAdmissibleModulesSymmetric1989} that it produces objects in the categorical center of the finite Hecke category attached to a reductive group $G$ and its Borel subgroup $B$. This observation has been elaborated upon in various contexts in (among other works) \cite{bezrukavnikovCharacterDmodulesDrinfeld2014}, \cite{ben-zviCharacterTheoryComplex2009}, \cite{lusztigTruncatedConvolutionCharacter2014}, \cite{lusztigNonunipotentCharacterSheaves2015}, \cite{holicharactersheveshomfly2023}, \cite{bitvcategoricalcenter2023}.

Another remarkable property of the Harish-Chandra transform on the character sheaves, playing a key role in \cite{bezrukavnikovCharacterDmodulesDrinfeld2014}, is that it becomes $t$-exact and Verdier self-dual after composing it with the Radon transform functor $R_*$. The latter transform is given by the averaging with respect to the unipotent radical $U^-$ of the opposite Borel subgroup $B^-$, or, equivalently, a pull-push functor along the following diagram:
\[
	(U\backslash G/U)/_{\Ad} T \xleftarrow{} (G/U)/_{\Ad} T \xrightarrow{} (U^{-}\backslash G/U)/_{\Ad} T,
\]
\[
	R_*\colon D^b \left((U\backslash G/U)/_{\Ad} T \right) \to D^b \left((U^-\backslash G/U)/_{\Ad} T \right).
\]
Here the subscript $*$ in $R_*$ indicates that the functor is defined using $*$-averaging.

For a local system $\mathcal{L}$ of finite order, let $D^b_{\mathfrak{C}_\mathcal{L}}(G)$ be the derived category of character sheaves with the character $\mathcal{L}$ (see Section \ref{sec:section-character-sheaves} for the precise definition). Our main result is as follows:
\begin{theorem}[Corollary \ref{RHC} and Remark \ref{main-remark}]\label{intro-main-result}
	The functor $R_* \circ \hc_!$, when restricted to the category $D^b_{\mathfrak{C}_\mathcal{L}}(G)$, is $t$-exact with respect to the perverse t-structure and commutes with the Verdier duality.
\end{theorem}
We prove Theorem \ref{intro-main-result} by identifying $R_*\circ\hc_!$ with the nearby cycles functor along the Vinberg degeneration. The key idea is to use that the nearby cycles functor commutes with the hyperbolic localization, as shown in \cite{nakajimalectures2016}, \cite{richarzspaceswithaction2019}. 

All our methods are geometric and can be restated with the same proofs, up to standard modifications, when $D^b(\cX)$ for stacks $\cX$ over $\mathbb{C}$ is replaced by either of the following:
\begin{itemize}
	\item the bounded derived category of holonomic $D$-modules,
	\item the bounded constructible derived category of $\Qlbar$-sheaves on a stack $\cX$ defined over an algebraically closed field of characteristic $p > 0$ for primes $p \neq \ell$ as in \cite{laszloSixOperationsSheaves2008a}. 
\end{itemize}
The parabolic case is treated in Section~\ref{sec:gener-pair-arbitr} by the same methods.

\subsection{Relation to other works.} The idea of relating the nearby cycles functor along the Vinberg degeneration to the functor $R_* \circ \hc_!$ already appears in \cite{bezrukavnikovCharacterDmodulesDrinfeld2014}, where an analogue of Theorem \ref{intro-main-result} in the setting of D-modules is proved. The same idea was used in \cite{chenFormulaGeometricJacquet2017} to prove Theorem \ref{intro-main-result} for unipotent character sheaves (that is, for the trivial local system $\mathcal{L}$), with methods applicable to either of the sheaf-theoretic cases, and also in the more general setting associated with a symmetric pair.  

A systematic treatment of when the hyperbolic localization commutes with various pull-push functors is given in \cite{linchennearbycycles2023}. We believe that an analogue of our Theorem~\ref{sec:character-sheaves-1} can be obtained using techniques developed in \emph{loc. cit}.

The fact that the composition of the Harish-Chandra and Radon transforms commutes with Verdier duality was also recently proved in \cite{grv26} by a completely different method (and in a different generality).

It seems likely that the methods of \cite{chenFormulaGeometricJacquet2017} can be adapted to treat the general monodromy $\mathcal{L}$ in the Borel case. However, we believe that the systematic use of hyperbolic localization presented in this paper gives a particularly short and straightforward proof. 

\subsection{Acknowledgments.} We would like to thank Roman Bezrukavnikov, Alexander Yom Din, and Yakov Varshavsky for helpful discussions. The work of R.G. was supported by the EPSRC grant EP/V053787/1. R.G. also acknowledges the hospitality of the Hebrew University of Jerusalem during his postdoctoral fellowship.  While working on the project, K. T. was at different times funded by the EPSRC program grant EP/R034826/1 and by the Deutsche Forschungsgemeinschaft SFB 1624 grant, Projektnummer 506632645.

\section{Setup and the statement of the main results}
\subsection{Limit functors.}\label{paragraph: limit functor} For a scheme $\tilX$ and a map $f \colon \tilX \to \mathbb{A}^{1}$, let 
\[
	\psi_f \colon D^b(\tilX \times_{\mathbb{A}^1}\mathbb{G}_m) \to D^b(\tilX \times_{\mathbb{A}^1} \operatorname{pt})
\] 
denote the nearby cycles functor. Here $\times_{\mathbb{A}^1}$ denotes the fiber product with respect to the embedding of the origin and the complementary open embedding $\operatorname{pt} \to \mathbb{A}^1 \leftarrow \mathbb{G}_m$, respectively. 

Let $X$ be a scheme with a $\mathbb{G}_m$-action, and denote the corresponding action map by $a \colon \mathbb{G}_m \times X\to X$. For $\tilX = X \times \mathbb{A}^{1}$, let $f \colon \tilX \to \mathbb{A}^1$ be the projection to the second factor. Following \cite{chenFormulaGeometricJacquet2017}, define the \emph{limit functor}
\[
  \Lim \colon D^b(X) \to D^b(X), \quad \Lim = \psi_f\circ a^*[1].
\]

Let $H$ be an algebraic group acting on $X$, and let us fix a cocharacter $\gamma \colon \mathbb{G}_m \to H$. This defines the $\mathbb{G}_m$-action on $X$, and we denote the corresponding limit functor by
\[
	\Limg \colon D^b(X) \to D^b(X).
\] 

Let $G$ be a reductive group and $T$ its maximal torus. A choice of a regular cocharacter $\lambda \colon \mathbb{G}_m \to T$ fixes a choice of a pair of opposite Borel subgroups $B, B^{-} \supset T$ for which $\lambda$ is dominant and antidominant, respectively. Let $U$ and $U^{-}$ be the unipotent radicals of $B$ and $B^{-}$, respectively. We have two projections $q_B \colon B \to T$ and $q_{B^{-}}\colon B^{-} \to T$.  

%and let $w_0 \in W$ be its longest element. 

Let $K \simeq G$ denote the diagonal subgroup of $G \times G$, and define three subgroups $K_0$, $K_\infty$, $K' \subset G \times G$ as the following fiber products
\begin{equation} \label{eq: def K}
	K_0 = B \times_T B^{-}, \quad K_\infty = B^{-}\times_T B, \quad K' = B \times_T B,
\end{equation}
with respect to the homomorphisms $q_{B}\colon B \to T$ and $q_{B^{-}}\colon B^- \to T$.

Fix a cocharacter of $H = G \times G$ given by
\[\gamma\colon \mathbb{G}_m \to G \times G, \quad \gamma = (\lambda^{-1}, \lambda).\]

\iffalse
\begin{remark}
 	Note that $B \times B^-$ and $B^- \times B$ are the repelling and attracting loci for the $\mathbb{G}_m$-action on $G\times G$ obtained by restricting the adjoint $G \times G$-action on $G \times G$ to $\mathbb{G}_m$ via $\gamma$. The subgroups $K_0 \subset B \times B^-$ and $K_{\infty} \subset B^- \times B$ are about as the limit of $K$ under the $\mathbb{G}_m$-action.
\end{remark}
\fi

Let $X$ be a scheme equipped with a left action of $G \times G$. As shown in \cite[Lemma 4.12]{chenFormulaGeometricJacquet2017}, the limit functor $\Limg$ associated with this case can be lifted to either of the following functors:
\begin{align}
	\Limg \colon D^b(K \backslash X) \to D^b(K_0 \backslash X), \label{eq: def L K0X}\\
	\Limg \colon D^b(K'\backslash X) \to D^b(K' \backslash X). \label{eq: def L K'X}
\end{align}

\begin{remark}
	Consider the constant group scheme $\widetilde{G} = G \times \mathbb{P}^1$ over $\mathbb{P}^1$. One can define its subgroup scheme $\widetilde{K}$, with the fibers $K_t = \gamma(t) K \gamma(t)^{-1}$ for $t \neq 0, \infty$, and the fibers over $0$ and $\infty$ given by $K_0$ and $K_{\infty}$; see \eqref{eq: def K}. One can verify that $\widetilde{K}$ is a well-defined smooth group scheme over $\mathbb{P}^1$, see \cite[Proposition 2.5.2]{drinfeldTheoremBraden2014}.

	Let $\widetilde{X} = X \times \mathbb{P}^1$ and $\widetilde{K}'= K' \times \mathbb{P}^1$ be the constant scheme/constant group scheme over $\mathbb{P}^1$. Informally, the functors \eqref{eq: def L K0X} and \eqref{eq: def L K'X} can be thought of as the limit functor for the families $\widetilde{K} \backslash \widetilde{X},$ $\widetilde{K}' \backslash \widetilde{X} \rightarrow \mathbb{P}^1$. We refer to \cite{chenFormulaGeometricJacquet2017} for details.
\end{remark}

We collect the required properties of the limit functor in the following lemmas. They follow from the standard properties of the nearby cycles functor, see \cite[XIII]{deligne2006groupes}.
\begin{lemma}\label{lemma:limit-functors1}
	Let $X_1, X_2$ be schemes equipped with a $\mathbb{G}_m$-action, and let $f \colon X_1 \to X_2$ be a $\mathbb{G}_m$-equivariant map.
	\begin{enumerate}[label=(\alph*),ref=(\alph*)]
		\item \label{item:1}   There are natural transformations of functors $D^b(X_1) \to D^b(X_2)$
		\[
			f_!\circ \Lim \to \Lim\circ f_!, \quad \quad
			\Lim\circ f_* \to f_* \circ \Lim.
		\]
		If $f$ is proper, both natural transformations are isomorphisms.
		\item\label{item:2} There are natural transformations of functors $D^b(X_2) \to D^b(X_1)$:
		\[
			\Lim \circ f^! \to f^!\circ\Lim, \qquad
			f^*\circ\Lim\to \Lim \circ f^*.
		\] 
  	If $f$ is smooth, both natural transformations are isomorphisms.
  	\item\label{item:3} The following diagram is commutative
		\begin{equation*}
		\begin{tikzcd}
			\Lim  \arrow[r] \arrow[d]  & \Lim \circ f^! f_! \arrow[d] \\
			f^! f_! \circ\Lim \arrow[r]  & f^! \circ \Lim \circ f_!         
		\end{tikzcd}
		\end{equation*}
	\end{enumerate}
\end{lemma}
 The next lemma uses the notion of $\mathbb{G}_m$-monodromic sheaves; we give a precise definition in Subsection \ref{sec:section-character-sheaves}.
\begin{lemma}\label{lemma:limit-functors2}
	When restricted to the subcategory of $\mathbb{G}_m$-monodromic sheaves on $X_1$, the functor $\Lim$ is isomorphic to the identity functor.
\end{lemma}
 Consider a Cartesian diagram 
\begin{equation} \label{eq: diag Cart}
	\begin{tikzcd}
		X_1 \times_S X_2 \arrow[r, "p_2"] \arrow[d, "p_1"'] & X_1 \arrow[d, "q_1"] \\
		X_2 \arrow[r, "q_2"']                               & S                   
	\end{tikzcd}
\end{equation}
of schemes with a $\mathbb{G}_m$-action, where all the arrows are $\mathbb{G}_m$-equivariant. 

\begin{lemma} \label{lemma:limit-functors3}
	The following diagram, where the horizontal arrows are induced by the natural transformations from Lemma~\ref{lemma:limit-functors1} and the vertical ones are base change isomorphisms, is commutative
  % https://tikzcd.yichuanshen.de/#N4Igdg9gJgpgziAXAbVABwnAlgFyxMJZABgBpiBdUkANwEMAbAVxiRAEcB9AJgD0AqLsACMAQgC+AHUkBjLACcZ0gLZ0cACwBGm4ABlxIcaXSZc+QimHkqtRizZc+-aXMUq1WnfpcKZAAiExAyMTbDwCIm5ranpmVkQQdw1tPXE-H0UAngEskQlDYxAMMPMiMmEbWPsEtE5gbglakXEBDKVJVWSvYMLiswiUKIqYu3jEjo8U-XTZXz8mhrSm4Rb+AtD+i2QrYds4tgWJNqTPVLb5upWBQxsYKABzeCJQADN5CGUkMhAcCCQVwpvD7-ai-JDcEIgIGfRDfMGIADMkOh4NBf0QABZke8YQi0UgAKzY4GIAn4zHiCjiIA
  \begin{equation}\label{eq:4}
  \begin{tikzcd}
		q_2^*q_{1!}\circ\Lim \arrow[r] \arrow[d] & q_2^*\circ\Lim\circ q_{1!} \arrow[r]   & \Lim \circ q_2^* q_{1!} \arrow[d] \\
		p_{1!}p_{2}^*\circ\Lim \arrow[r]         & p_{1!}\circ\Lim\circ p_{2}^* \arrow[r] & \Lim \circ p_{1!} p_{2}^*        
	\end{tikzcd}
	\end{equation}
\end{lemma}
\begin{remark}
	It may be convenient for further arguments to think about the commutativity of diagram \eqref{eq:4} as follows. Consider the cube diagram below. Its front and back faces are identical and induced by the Cartesian square \eqref{eq: diag Cart}. Its arrows coming out of the plane of the paper are the limit functors $\Lim$. 
  %% Adapted from https://tex.stackexchange.com/questions/256835/tikz-cd-shade-faces-of-commutative-cube
	\[
		\definecolor{cof}{RGB}{219,144,71}
		\definecolor{pur}{RGB}{186,146,162}
		\definecolor{greeo}{RGB}{52,111,72}
		\definecolor{greet}{RGB}{52,111,72}
		\definecolor{red-top}{RGB}{186,146,162}
		\definecolor{red-back}{RGB}{140,120,100}
		\definecolor{red-right}{RGB}{160,133,140}
		\definecolor{green-front}{RGB}{91,173,69}
		\definecolor{green-left}{RGB}{52,111,72}
		\definecolor{green-bottom}{RGB}{73,140,70}
		\tikzset{FSTYLE/.style = {draw=blue, % just to demonstrate, where LA is used
		        line width=#1, -{Straight Barb[length=3pt]}},
		        FSTYLE/.default=1pt
		        }
		\tikzset{BSTYLE/.style = {draw=red, % just to demonstrate, where LA is used
		        line width=#1, -{Straight Barb[length=3pt]}},
		        BSTYLE/.default=1pt
		        }
		\begin{tikzcd}
			% [execute at end picture={
			%   \foreach \Nombre in  {A,B,C,D,E,F,G,H,I,J}
			%   {\coordinate (\Nombre) at (\Nombre.center);}
			%   \fill[red-back,opacity=0.2] 
			%   (B) -- (C) -- (G) -- (F) -- cycle;
			%   \fill[green-left,opacity=0.2] 
			%   (A) -- (B) -- (F) -- (E) -- cycle;
			%   \fill[green-bottom,opacity=0.2] 
			%   (E) -- (F) -- (G) -- (H) -- cycle;
			%   \fill[red-top,opacity=0.2] 
			%   (A) -- (B) -- (C) -- (D) -- cycle;
			%   \fill[green-front,opacity=0.2] 
			%   (A) -- (D) -- (H) -- (E) -- cycle;
			%   \fill[red-right,opacity=0.2] 
			%   (H) -- (D) -- (C) -- (G) -- cycle;
			% }]
			&
			|[alias=B]|D^b(X_1)
			\ar[FSTYLE, dashed]{dl}[swap]{\Lim}
			\ar{rr}{q_{1!}}
			\ar[BSTYLE]{dd}[near end]{}
			& & |[alias=C]|D^b(S)
			\ar{dd}{q_2^*}
			\ar{dl}[swap, sloped, near start]{}
			\\
			|[alias=A]|D^b(X_1)
			\ar[crossing over, FSTYLE, dashed]{rr}[near start]{}
			\ar{dd}[swap]{p_2^*}
			& |[alias=I]| & |[alias=D]|D^b(S)
			\\
			&
			|[alias=F]|D^b(X_1\times_S X_2)
			\ar[BSTYLE, near start]{rr}{}
			\ar[sloped, near end]{dl}{}
			& |[alias=J]| & |[alias=G]|D^b(X_2)
			\ar[BSTYLE]{dl}{\Lim}
			\\
			|[alias=E]|D^b(X_1\times_S X_2)
			\ar{rr}[swap]{p_{1!}}
			& & |[alias=H]|D^b(X_2)
			\ar[crossing over, FSTYLE, from=uu, dashed]{}{}
		\end{tikzcd}
	\]
	Lemma \ref{lemma:limit-functors3} says that two natural transformations between red and dashed blue paths $D^b(X_1) \to D^b(X_2)$, one given by composition of natural transformations on the top, right and back faces of the cube, and another along front, left, and bottom faces, coincide. 
\end{remark}

\subsection{Character sheaves.}\label{sec:section-character-sheaves}
Let $X$ be a scheme with an action of an algebraic torus $T$, and let $a \colon T \times X \to X$ be the action map. Let $\mathcal{L}$ be a rank 1 local system on $T$. We say that $\F \in D^b(X)$ is \emph{weakly $\mathcal{L}$-equivariant} if there is an isomorphism $a^*\F \simeq \mathcal{L} \boxtimes \F$. We denote by $D^b_{\mathcal{L}}(X)$ the full idempotent-complete triangulated subcategory of $D^b(X)$ generated by weakly $\mathcal{L}$-equivariant complexes. Its objects will be called \emph{$\mathcal{L}$-monodromic} or \emph{monodromic complexes with generalized monodromy $\mathcal{L}$}. 

We call a local system $\mathcal{L}$ on $T$ \emph{rational} if its pullback $\eta^*\mathcal{L}$ is trivial for some surjective homomorphism $\eta \colon T \to T$. 

If $T = \mathbb{G}_m$, we will denote by $D^b_{mon}(X)$ the category generated by $D^b_{\mathcal{L}}(X)$ for all rational $\mathcal L$ and refer to this category as the category of \emph{$\mathbb{G}_m$-monodromic} complexes.

For a scheme $X$ with an action of an algebraic group $H$ which is normalized by the action of $T$, we will write $D^b_{\mathcal{L}}(H\backslash X)$ for the full triangulated subcategory of $D^b(H \backslash X)$ consisting of objects $\F$ with $\For^H\F \in D^b_{\mathcal{L}}(X)$.

Let $X = G \times G$, and consider the actions of $H = G\times G$ on $X$ by left and right diagonal multiplication:
\[(g_1, g_2) \big((x_1, x_2)\big) = (g_1 x_1 g_2^{-1}, g_1 x_2 g_2^{-1}).\]
Notice that $H \backslash X \simeq K \backslash X / K$. Consider a functor
\[
	\hc_! \colon D^b(K\backslash X/K) \to D^b(K\backslash X/K'), \quad \hc_! = \prescript{R}{}{\Av}^{K'}_{K \cap K',!}\circ\prescript{R}{}{\For}^{K}_{K\cap K'},
\]
where the superscript $R$ on averaging and forgetful functors is used to emphasize that these are with respect to the groups acting on the right. The notation $\hc$ stands for Harish-Chandra, and the functor thus defined is commonly called the \emph{Harish-Chandra functor}.  
Its right adjoint functor is, up to a homological shift,
\[
	\chi \colon D^b(K\backslash X/K') \to D^b(K\backslash X/K), \quad \chi = \prescript{R}{}{\Av}^{K}_{K \cap K',!}\circ\prescript{R}{}{\For}^{K'}_{K\cap K'}.
\]
Indeed, $\prescript{R}{}{\Av}^{K}_{K \cap K',!}$ is a proper pushforward, so $*$ and $!$-versions coincide.
\begin{remark}
	For a smooth map the $!$-pullback and $*$-pullback differ by the standard homological shift. To keep the notation light, we suppress this shift and write simply $\For$. The reader who wishes to keep track of the precise grading can recover the shift from whether the corresponding forgetful functor is realized as $!$- or $*$-pullback.
\end{remark}

The torus $T$ acts on $K\backslash X/K'$ on the right. Let $\mathcal{L}$ be a rational rank 1 local system on $T$.
Recall that a strictly full triangulated subcategory is called \emph{thick} if it is closed under direct summands. The derived category of \emph{character sheaves} with character $\mathcal{L}$, denoted by $D^b_{\CS_{\mathcal{L}}}(G)$, is the thick subcategory of $D^b(K\backslash X/K)$ generated by the essential image of $D^b_{\mathcal{L}}(K \backslash X/K')$ under the functor $\chi$.

Let $W = N_G(T)/T$ be the Weyl group of $G$. It is well known that the subcategories $D^b_{\CS_{\mathcal{L}}}(G), D^b_{\CS_{\mathcal{L}'}}(G)$ coincide if $\mathcal{L}$ and $\mathcal{L}'$ are in the same $W$-orbit, and otherwise are orthogonal. We denote by $D^b_{\CS}(G)$ the full triangulated subcategory of $D^b(K\backslash X/K)$ generated by $D^b_{\CS_{\mathcal{L}}}(G)$ for all rational $\mathcal{L}$.

\subsection{Main results.}
Consider the diagram:
\begin{equation}\label{eq:1}
	\begin{tikzcd}
		D^b(K\backslash X/K) \arrow[r, "\Limg"]\arrow[d, "\Av^{K'}_{B,!}\circ \For^K_B"'] & D^b(K_0\backslash X/K) \arrow[d,"\Av^{K'}_{B,!}\circ \For^{K_0}_B"] \arrow[dl, Rightarrow, "\alpha"]\\
		D^b(K'\backslash X/K) \arrow[r, "\Limg"'] & D^b(K'\backslash X/K)
   \end{tikzcd}
\end{equation}

Here $\alpha \colon \Av^{K'}_{B,!}\circ \For^{K_0}_B\circ\Limg \to \Limg\circ\Av^{K'}_{B,!}\circ \For^K_B$ is the natural transformation given by Lemma \ref{lemma:limit-functors1}. We will prove the following.
\begin{theorem}\label{sec:character-sheaves}
	The restriction of $\alpha$ to the category of character sheaves $D^b_{\CS}(G)$ is an isomorphism.
\end{theorem}
\begin{corollary}\label{RHC}
	After restriction to the category of character sheaves $D^b_{\CS}(G)$, there is an isomorphism of functors
	\[
		\Limg \simeq \Av^{K_0}_{B,*}\circ \For^{K'}_B \circ \Av^{K'}_{B,!}\circ \For^K_B.
	\]
\end{corollary}
\begin{proof}
	By Lemma \ref{lemma:limit-functors2}, the lower horizontal arrow in the diagram \eqref{eq:1} is isomorphic to the identity functor when restricted to the monodromic subcategory. We obtain an isomorphism
	\[
		\Av^{K'}_{B,!}\circ \For^{K_0}_B \circ  \Limg \simeq \Av^{K'}_{B,!}\circ \For^K_B.
	\]
	It is well known that, restricted to the monodromic subcategory, the functor $\Av^{K'}_{B,!}\circ \For^{K_0}_B$ is an equivalence, with its inverse given by $\Av^{K_0}_{B,*}\circ \For^{K'}_B$, and the corollary follows. 
\end{proof}
\begin{remark}\label{main-remark} 
	To match the notation of the introduction, note that 
	\[
		K\backslash X/K = G\backslash G \times G/G \simeq G/_{\Ad}G, 
	\]
  	\[ 
  		K'\backslash X/K = (B \times_T B)\backslash G \times G/G \simeq (U\backslash G/U)/_{\Ad}T, 
	\]
	\[
		K_0\backslash X/K = (B \times_T B^-)\backslash G \times G/G \simeq (U^-\backslash G/U)/_{\Ad}T.
	\]
\end{remark} 

\bigskip

We will reduce Theorem \ref{sec:character-sheaves} to a similar statement of the form ``the nearby cycles functor commutes with the averaging functor'' with $X/K$ replaced first by $X/K'$ and then by $\fl \times \fl$, where $\fl$ is the flag variety of $G$. 
Consider the following diagram, which is similar to diagram \eqref{eq:1}, with $X/K$ replaced by $X/K'$.
\begin{equation}\label{eq:2}
	\begin{tikzcd}
		D^b(K\backslash X/K') \arrow[r, "\Limg"]\arrow[d, "\Av^{K'}_{B,!}\circ \For^K_B"'] & D^b(K_0\backslash X/K') \arrow[d,"\Av^{K'}_{B,!}\circ \For^{K_0}_B"] \arrow[dl, Rightarrow, "\alpha'"]\\
		D^b(K'\backslash X/K') \arrow[r, "\Limg"'] & D^b(K'\backslash X/K')
	\end{tikzcd}
\end{equation}
Here $\alpha' \colon \Av^{K'}_{B,!}\circ \For^{K_0}_B\circ\Limg \to \Limg\circ\Av^{K'}_{B,!}\circ \For^K_B$ is again the natural transformation given by Lemma \ref{lemma:limit-functors1}. We have the following.
\begin{theorem}\label{sec:character-sheaves-1}
	For any rational rank-one local system $\mathcal{L}$ on $T$, the restriction of $\alpha'$ to the monodromic category  $D^b_{\mathcal{L}}(K\backslash X/K')$ is an isomorphism. 
\end{theorem}

\bigskip

Let $\fl = G/B$ be the flag variety of $G$. Consider the following diagram, which is similar to diagram \eqref{eq:1}, with $X/K'$ replaced by $\fl \times \fl$.
\begin{equation}\label{eq:3}
	\begin{tikzcd}
		D^b(K\backslash \fl\times\fl) \arrow[r, "\Limg"]\arrow[d, "\Av^{K'}_{B,!}\circ \For^K_B"'] & D^b(K_0\backslash \fl\times\fl) \arrow[d,"\Av^{K'}_{B,!}\circ \For^{K_0}_B"] \arrow[dl, Rightarrow, "\alpha''"]\\
		D^b(K'\backslash \fl\times\fl) \arrow[r, "\Limg"'] & D^b(K'\backslash \fl\times\fl)
	\end{tikzcd}
\end{equation}
Here $\alpha'' \colon \Av^{K'}_{B,!}\circ \For^{K_0}_B\circ\Limg \to \Limg \circ\Av^{K'}_{B,!}\circ \For^K_B$ is again the natural transformation given by Lemma \ref{lemma:limit-functors1}.
\begin{theorem}\label{sec:character-sheaves-2}
	The transformation $\alpha''$ is an isomorphism. 
\end{theorem}

The implication Theorem \ref{sec:character-sheaves-1} $\Rightarrow$ Theorem \ref{sec:character-sheaves} will be deduced from a diagram chase using the standard properties of the nearby cycles functor. For unipotent monodromy, Theorem \ref{sec:character-sheaves-2} implies Theorem \ref{sec:character-sheaves-1} in a similar way. The proof of Theorem~\ref{sec:character-sheaves-1} and Theorem~\ref{sec:character-sheaves-2} will use properties of the hyperbolic localization functor, which we now recall.

\section{Hyperbolic localization}
\subsection{Braden's theorem.} We recall some notions from \cite{drinfeldTheoremBraden2014}. Let $X$ be a scheme and $a \colon \mathbb{G}_m \times X \rightarrow X$ be a $\mathbb{G}_m$-action. Consider the fixed points scheme $X^{\mathbb{G}_m}$. Let $X^+$ and $X^-$ be the attracting and the repelling schemes, respectively. Consider the diagram
\[
	\begin{tikzcd}
		& \arrow[dl, "{p_{\pm}}"'] X^{\pm} \arrow{dr}{i_{\pm}} \\
		X^{\mathbb{G}_m}   && X
	\end{tikzcd}
\]
where $i_{\pm}$ is the inclusion and $p_{\pm}$ is the limit map. Let us define the functors of hyperbolic restriction
\begin{align}
	\Hf^{+} =& (p_{+})_! \circ i_+^* & \Hf^{-} =& (p_{-})_* \circ i_-^!
\end{align}

\begin{theorem}[\cite{BradenHyperbolic03}, \cite{drinfeldTheoremBraden2014}]
	There is a natural transformation $\Hf^{-} \rightarrow \Hf^{+}$ which is an isomorphism on $\mathbb{G}_m$-monodromic complexes.
\end{theorem}
Denote by $\Hf$ the restriction of $\Hf^{\pm}$ to the $\mathbb{G}_m$-monodromic subcategory.

\subsection{Hyperbolic localization commutes with nearby cycles.} For a scheme $\tilX$ and a map $f \colon \tilX \to \mathbb{A}^{1}$ one has the nearby cycles functor $\psi_f$ as in Subsection \ref{paragraph: limit functor}. Let $a' \colon \mathbb{G}_m \times \tilX \rightarrow \tilX$ be an action which commutes with the map $f$ for the trivial action on $\mathbb{A}^1$.  Consider the fixed points of the action and the map $f_0 \colon \tilde{X}^{\mathbb{G}_m} \rightarrow \mathbb{A}^1$. The natural transformations from Lemma \ref{lemma:limit-functors1} yield natural transformations
\begin{align}\label{eq: Hpm and nearby}
	&\Hf^+ \circ \psi_f \rightarrow \psi_{f_0} \circ \Hf^+ &
	&\psi_{f_0} \circ \Hf^- \rightarrow \Hf^- \circ \psi_f
\end{align}

\begin{theorem}\label{thm: commutes with nearby cycles}
	The natural transformations~\eqref{eq: Hpm and nearby} yield an
	isomorphism
	\[
		\Hf\circ\psi_f\simeq\psi_{f_0}\circ\Hf
	\]
	upon restriction to the category
	of $\mathbb G_m$-monodromic sheaves
	$D^b_{mon}(\tilde X\times_{\mathbb A^1}\mathbb G_m)$.
\end{theorem}

\begin{proof}
	Let us consider $\mathbb G_m$-monodromic sheaves with the generalized monodromy $\mathcal L$. The case of the trivial $\mathcal L$ was established in \cite[Theorem 3.3]{richarzspaceswithaction2019}; see also \cite[Proposition~5.4.1(2)]{nakajimalectures2016}. For non-trivial $\mathcal L$ the theorem is evident since the hyperbolic restriction equals zero for both sides.
\end{proof}

%Informally, one can refer to the theorem above as \emph{hyperbolic restriction commutes with nearby cycles}.

\subsection{Hyperbolic localization commutes with Radon transform. Case of the flag variety.}
Consider the $\mathbb{G}_m$-action on $\fl\times\fl$ given by a cocharacter $(\lambda, \lambda)$ for a dominant $\lambda$. This gives us the functor
\[
	\Hf \colon D^b_{mon}(\fl \times \fl) \rightarrow  D^b(\fl^{\mathbb{G}_m} \times \fl^{\mathbb{G}_m})
\]
where the subscript $mon$ indicates the $\mathbb G_m$-monodromic subcategory for the $\mathbb G_m$-action defined above. Moreover, consider $U$ acting on the right factor of $\fl\times\fl$. This gives us the functors
\begin{align*}
	\For^{U} \colon&D^b_{mon}( \fl \times U \backslash \fl) \rightarrow D^b_{mon}( \fl \times \fl)\\
	\Av_!^{U} \colon&D^b_{mon}(\fl \times \fl) \rightarrow D^b_{mon}(  \fl \times  U \backslash \fl) 
\end{align*}
The following lemma is straightforward; see \cite[Sect. 1.4.5]{chenCasselmanJacquetFunctor2019}.
\begin{lemma}\label{lemma: CGY}
	The unit map induces a natural isomorphism.
	$$ \Hf \simeq \Hf \circ \For^{U} \circ \Av_!^{U}$$
\end{lemma}
As a corollary we obtain the following proposition.
\begin{proposition}\label{prop: radon commutes with H}
	\begin{enumerate}[label=(\alph*),ref=(\alph*)]
		\item\label{item: K} The unit map $\Id \rightarrow \For^{K'}_B  \circ \Av^{K'}_{B,!}$ induces an isomorphism of functors making the following diagram commutative:
		\begin{equation*}
		\begin{tikzcd}
			D^b (K \backslash \fl\times\fl) \arrow[rd, "\Hf \circ \For^{K}"] \arrow[rr, "\Av^{K'}_{B,!} \circ \For^{K}_B"] && \arrow[ld, "\Hf \circ \For^{K'}"] D^b(K'\backslash \fl\times\fl)  \\ 
			& D^b(\fl^{\mathbb{G}_m} \times\fl^{\mathbb{G}_m})
		\end{tikzcd}
		\end{equation*}
		\item\label{item: K_0} The same holds for $K$ replaced by $K_0$.
	\end{enumerate}
\end{proposition}
\begin{proof}
	First of all, notice that the sheaves in the essential image of $\For^{K}$ and $\For^{K'}$ are $\mathbb{G}_m$-monodromic. We used this implicitly in the formulation of the proposition since we wrote $\Hf$ rather than~$\Hf^{\pm}$.
	
	Using the base change isomorphism for the following Cartesian diagram
	\begin{equation}\label{base_change:flags}
	\begin{tikzcd}
		\fl \times \fl \arrow[r] \arrow[d] & B \backslash \fl \times \fl \arrow[d] \\
		\fl \times U \backslash \fl \arrow[r] & K' \backslash \fl \times \fl                  
	\end{tikzcd}
	\end{equation}
	we get
	\begin{multline*}
		\For^{K'} \circ \Av^{K'}_{B,!} \circ \For^{K}_B \simeq 
		\For^{U} \circ \For^{K'}_{U} \circ \Av^{K'}_{B,!} \circ \For^{K}_B \simeq \\
		\For^{U} \circ \Av^{U}_! \circ \For^B \circ \For^{K}_B \simeq \For^{U} \circ \Av^{U}_! \circ \For^{K}.
	\end{multline*}
	Again, since the sheaves in the essential image of $\For^{K}$ are $\mathbb{G}_m$-monodromic, we can apply Lemma \ref{lemma: CGY}:
	\begin{equation}\label{eq:nat_transform:flags}
		\Hf \circ \For^{K} \xrightarrow{\sim} 
		\Hf \circ \For^{U} \circ \Av^{U}_! \circ \For^K \simeq 
		\Hf \circ \For^{K'} \circ \Av^{K'}_{B,!} \circ \For^{K}_B.
	\end{equation}
	It remains to show that \eqref{eq:nat_transform:flags} is indeed induced by the unit map $\Id \rightarrow \For^{K'}_B  \circ \Av^{K'}_{B,!}$. More precisely, the latter natural transformation is given by
	\begin{multline*}
		\Hf \circ \For^K \simeq \Hf \circ \For^B \circ \For^K_B \xrightarrow{\sim} \\  
		\Hf \circ \For^B \circ \For^{K'}_B  \circ \Av^{K'}_{B,!} \circ \For^K_B \simeq
		\Hf \circ \For^{K'}  \circ \Av^{K'}_{B,!} \circ \For^K_B.
	\end{multline*}
	These natural transformations are equal since smooth base-change commutes with the unit transformation; see diagram \eqref{base_change:flags}. Hence part~\ref{item: K} follows. The case of $K_0$ is similar.
\end{proof}
%Informally, one can refer to the proposition above as \emph{hyperbolic restriction commutes with Radon transform $R_! = \Av^{K'}_{B,!} \circ \For^{K_0}_B$}.

\subsection{Hyperbolic localization commutes with Radon transform. Case of the base affine space.}
For $w_1,w_2 \in W$, consider the $\mathbb{G}_m$-action on $X/K'$ given by
\begin{gather}
	\aw\colon \mathbb{G}_m \times X/K' \rightarrow X/K'\\ 
	\aw\left(t, (g_1, g_2) \right) = \Big( \lambda(t) \cdot g_1 \cdot w_1^{-1} \left( \lambda (t^{-1}) \right)  ,
	\lambda(t) \cdot g_2 \cdot w_2^{-1} \left(\lambda \left(t^{-1} \right) \right)\Big) \label{eq:w1w2_action}
\end{gather}
where $\lambda$ is a regular dominant cocharacter. This gives us the hyperbolic localization functors
\[
	\Hf_{w_1, w_2} \colon D^b_{mon}(X/K') \rightarrow  D^b((X/K')^{\mathbb{G}_m}),
\]
where the subscript $mon$ indicates the $\mathbb G_m$-monodromic subcategory for the $\mathbb G_m$-action $\aw$. The following lemma follows directly from the Bruhat decomposition.
\begin{lemma} For the $\mathbb{G}_m$-action on $X/K'$ given by $\aw$, the following hold.
	\begin{enumerate}[label=(\alph*),ref=(\alph*)]
	\item The fixed locus
	\[
	(X/K')^{\mathbb G_m} =
	\bigsqcup_{\tilde w_1,\tilde w_2} (\tilde w_1B\times\tilde w_2B)/K',
	\]
	where the disjoint union is taken over the pairs $\tilde{w}_1,\tilde{w}_2 \in W$ satisfying
	\begin{equation}\label{eq:involved_condition}
	w_1^{-1}(\lambda) - w_2^{-1}(\lambda) = \tilde{w}_1^{-1}(\lambda) - \tilde{w}_2^{-1}(\lambda).
	\end{equation}
	\item The attracting locus for $(\tilde{w}_1 B \times \tilde{w}_2 B)/K'$ is $(U\tilde{w}_1 B \times U\tilde{w}_2 B)/K'$.
	\end{enumerate}
\end{lemma}
For the rest of this subsection, we fix $w_1, w_2$ and shorten $\Hf := \Hf_{w_1,w_2}$. Consider $U$ acting on the right factor of $X/K' = \left(G \times G\right)/K'$. This gives us functors
\begin{align*}
	\For^{U}\colon&D^b_{mon}( U \backslash X/K') \rightarrow D^b_{mon}( X/K'), \\
	\Av_!^{U}\colon&D^b_{mon}(X/K') \rightarrow D^b_{mon}( U \backslash X/K').
\end{align*}
The following lemma is a straightforward generalization of \cite[Sect. 1.4.5]{chenCasselmanJacquetFunctor2019}.
\begin{lemma}\label{lemma: CGY monodromic} The unit map induces an isomorphism of functors
	\[\Hf \simeq \Hf \circ \For^{U} \circ \Av_!^{U}\]
\end{lemma}
\begin{proof}
	Recall that $p_+$ is the limit map for the attracting set. This map factors as
	\[
		\left(X/K'\right)^+ \rightarrow 
		U \backslash\left(X/K'\right)^+ \xrightarrow{\tilde{p}_+}
		\left(X/K'\right)^{\mathbb{G}_m} 
	\]
	Hence 
	\[
		(p_{+})_! \circ \For^{U} \circ \Av_!^{U} =  (\tilde{p}_{+})_!  \circ \Av_!^{U} \circ \For^{U} \circ \Av_!^{U} \simeq (\tilde{p}_{+})_! \circ \Av_!^{U} = (p_{+})_!.
	\]
	 Then we have the following chain of isomorphisms
	\begin{equation*}
		\Hf \circ \For^{U} \circ \Av_!^{U} = 
		(p_{+})_! \circ i_+^* \circ  \For^{U} \circ \Av_!^{U} \underset{(i)}{\simeq} 
		(p_{+})_! \circ \For^{U} \circ \Av_!^{U} \circ i_+^*   \underset{(ii)}{\simeq}
		(p_{+})_!  \circ i_+^*  = \Hf,
	\end{equation*}
	where $(i)$ is the base change isomorphism, and $(ii)$ is constructed above.
\end{proof}
As a corollary, we obtain the following proposition:
\begin{proposition} \label{prop: radon commutes with H monodromic}
	\begin{enumerate}[label=(\alph*),ref=(\alph*)]
		\item\label{item: K monodromic} The unit map $\Id \rightarrow \For^{K'}_B  \circ \Av^{K'}_{B,!}$ induces an isomorphism of functors making the following diagram commutative:
		\begin{equation*}
		\begin{tikzcd}
			D^b_{\mathcal L}(K \backslash X/K') \arrow[rd, "\Hf \circ \For^{K}"] \arrow[rr, "\Av^{K'}_{B,!} \circ \For^{K}_B"] && \arrow[ld, "\Hf \circ \For^{K'}"] D^b_{\mathcal L}(K'\backslash X/K')  \\ 
			& D^b((X/K')^{\mathbb{G}_m})
		\end{tikzcd}
		\end{equation*}
		where the subscript $\mathcal L$ indicates the $\mathcal L$-monodromic subcategory with respect to the right $T$-action.
		\item\label{item: K_0 monodromic} The same holds for $K$ replaced by $K_0$.
	\end{enumerate}
\end{proposition}

We use $\Hf$ in the proposition above because the essential images of the functors $\For^K$ and $\For^{K'}$ on the respective source categories lie in $D^b_{\mathrm{mon}}(X/K')$. The proof is the same as that of
Proposition~\ref{prop: radon commutes with H} with Lemma~\ref{lemma: CGY monodromic} in place of
Lemma~\ref{lemma: CGY}.

\iffalse
\begin{proof}
	Using base change for the following Cartesian diagram
	\[
	\begin{tikzcd}
		X/K' \arrow[r] \arrow[d]     & B \backslash X/K' \arrow[d] \\
		U \backslash X/K' \arrow[r]  & K' \backslash X/K'                 
	\end{tikzcd}
	\]
	we get %the following chain of isomorphisms of functors from $D^b(K \backslash X/K')$ to $D^b(X/K')$
	\begin{multline*}
		\For^{K'}\circ \Av^{K'}_{B,!} \circ \For^{K}_B \simeq 
		\For^{U} \circ \For^{K'}_{U} \circ \Av^{K'}_{B,!} \circ \For^{K}_B \simeq \\
		\For^{U} \circ \Av^{U}_! \circ \For^B \circ \For^{K}_B \simeq
		\For^{U} \circ \Av^{U}_! \circ \For^{K}
	\end{multline*}
	%The sheaves in the essential image of $\For^{K}$ are $\mathbb{G}_m$-monodromic on the attracting locus, which allows us to use 
	Composing the above isomorphism with $\Hf$ and applying Lemma~\ref{lemma: CGY monodromic} proves part~\ref{item: K monodromic}. The case of $K_0$ is identical. 
\end{proof}
\fi
In the proposition above, we considered the functor
$\Av^{K'}_{B,!}\circ \For^{K_0}_B$. This functor is closely related to
the Radon transform $R_!$; cf.
Subsection~\ref{subsec:intro_results}. These functors are given by the same formula, but in slightly different contexts. For this reason, we informally refer to $\Av^{K'}_{B,!}\circ \For^{K_0}_B$ as the Radon
transform, and to the proposition above as the statement that \emph{hyperbolic restriction commutes with the Radon transform}.
\begin{remark}
	We formulated Proposition \ref{prop: radon commutes with H monodromic} in the way we are going to use it. The proposition admits a straightforward generalization with $B$ taken instead of $K$ or $K_0$; this version implies both parts \ref{item: K monodromic} and \ref{item: K_0 monodromic}. Also, one does not have to restrict to the monodromic subcategories; the non-monodromic counterpart requires writing $\Hf^+$ instead of $\Hf$.
\end{remark}

\section{Proof of the main results}
\subsection{Proof for the flag variety.}\label{subsec:proof_flag}
In this subsection we give a proof of Theorem \ref{sec:character-sheaves-2}. The following lemma is straightforward to verify.
\begin{lemma}
The functor 
\[
	\Hf \circ \For^{K'}\colon D^b(K'\backslash \fl\times\fl) \rightarrow D^b(\fl^{\mathbb{G}_m} \times\fl^{\mathbb{G}_m})
\]
is conservative.
\end{lemma}
Using the lemma above, it is enough to check that 
\begin{equation} \label{eq: the natural transform}
	\Hf \circ \For^{K'}(\alpha'') \colon \Hf \circ \For^{K'} \circ \Av^{K'}_{B,!}\circ \For^{K_0}_B\circ\Limg \to  
	\Hf \circ \For^{K'} \circ \Limg \circ\Av^{K'}_{B,!}\circ \For^K_B
\end{equation}
is an isomorphism. We have the following chain of isomorphisms:
\begin{multline} \label{eq: funct iso1}
	\Hf \circ \For^{K_0} \circ \Limg \overset{(i)}{\simeq}
	\Hf \circ \Limg  \circ \For^{K} \overset{(ii)}{\simeq} \\
	\Limg \circ \Hf \circ \For^{K} \overset{(iii)}{\xrightarrow{\sim}} 
	\Limg  \circ \Hf \circ \For^{K'} \circ\Av^{K'}_{B,!}\circ \For^K_B \overset{(iv)}{\simeq} \\
	\Hf \circ \Limg \circ \For^{K'} \circ\Av^{K'}_{B,!}\circ \For^K_B \overset{(v)}{\simeq}
	\Hf \circ  \For^{K'} \circ \Limg  \circ\Av^{K'}_{B,!}\circ \For^K_B
\end{multline}
The first and the fifth isomorphisms are given by Lemma \ref{lemma:limit-functors1}, the second and the fourth by Theorem~\ref{thm: commutes with nearby cycles}. The third isomorphism is induced by the unit transformation, see Proposition~\ref{prop: radon commutes with H}\ref{item: K}.

Consider also the following isomorphism
\begin{equation}\label{eq: funct iso2}
	\Hf \circ \For^{K_0} \circ \Limg \xrightarrow{\sim}  
	\Hf \circ \For^{K'} \circ \Av^{K'}_{B,!}\circ \For^{K_0}_B\circ\Limg
\end{equation}
given by the unit transformation, see Proposition \ref{prop: radon commutes with H}\ref{item: K_0}. 

Lemma \ref{lemma:limit-functors1}\ref{item:3} implies that the natural transformation \eqref{eq: funct iso1} equals the composition of \eqref{eq: funct iso2} and \eqref{eq: the natural transform}. Since \eqref{eq: funct iso1} and \eqref{eq: funct iso2} are isomorphisms, so is \eqref{eq: the natural transform}, which completes the proof of Theorem \ref{sec:character-sheaves-2}.

\begin{remark}
	In the proof above, we used several functors which turned out to be equivalent. One can also notice that they are equivalent to $\Limg \circ \Hf \circ \For^{K} \simeq \Hf \circ \For^{K}$ since $\fl^{\mathbb{G}_m} \times \fl^{\mathbb{G}_m}$ is a trivial family.
\end{remark}

\begin{remark}
	Consider the diagonal embedding $\delta\colon K\backslash\fl \to K\backslash\fl \times \fl$. The Beilinson--Bernstein localization of the object $\mathbb{L}_\gamma(\delta_* \cc_{K\backslash\fl})$ is known as a semiregular module; see \cite{donin2008lie}. For $w \in W$, let $\cO_w = UwB \subset \fl$ and $\cO_w^{-} = U^{-}wB \subset \fl$ be the Bruhat cells, and let $i_w\colon\cO_w \to \fl$ and $i'_w\colon\cO^{-}_w \to \fl$ be the corresponding locally closed embeddings. Let $\Delta_w = i_{w!}\,\cc_{\cO_w}[\dim \cO_w]$ and $\nabla_w^- = i'_{w*}\,\cc_{\cO^{-}_w}[\dim \cO^-_{w}]$. Theorem \ref{thm: commutes with nearby cycles} implies that $\For^{K_0}\mathbb{L}_\gamma(\delta_* \cc_{K\backslash\fl})[\dim K_0]$ admits a filtration with associated graded pieces given by $\Delta_{w_0w^{-1}}\!\boxtimes \nabla^{-}_{w}$ for $w \in W$. After Beilinson--Bernstein localization, this recovers \cite[Theorem 11.14]{donin2008lie}. See also \cite{linchennearbycycles2023} for a related discussion in the affine setting.
\end{remark}

\subsection{Proof for monodromic sheaves.}
In this subsection we prove Theorem~\ref{sec:character-sheaves-1}. To do this, we use the hyperbolic localization $\Hf_{w_1, w_2}$ for the action \eqref{eq:w1w2_action}. Recall that the fixed locus for this action is the disjoint union of the components $(\tilde{w}_1B\times \tilde{w}_2B)/K'$, where $\tilde{w}_1$ and $\tilde{w}_2$ satisfy condition~\eqref{eq:involved_condition}. Consider the restriction of $\Hf_{w_1, w_2}$
\[
	\Hf^{\mathrm{res}}_{w_1, w_2}\colon D^b_{\mathcal L}(X/K') \to D^b \big( w_1B \times w_2B/K' \big).
\]
Denote $\Hfb = \bigoplus_{w_1, w_2 \in W} \Hf_{w_1, w_2}^{\mathrm{res}}$.
\begin{lemma}\label{sec:proof-monodr-sheav}
	The functor 
		\[
		\Hf_{\bullet} \circ \For^{K'}\colon D^b_{\mathcal L}(K'\backslash X/K') \rightarrow \bigoplus_{w_1, w_2} D^b \big( w_1B \times w_2B/K' \big)
		\]
		is conservative.
\end{lemma}
\begin{proof}
	For $(w_1,w_2) \in W^2$ and the corresponding action \eqref{eq:w1w2_action}, the attractor of the component $w_1B\times w_2B/K'$ is
	\[
	\cO_{w_1,w_2} = Uw_1B\times Uw_2B/K'.
	\]
	Denote the embedding and the contraction map by
	\[
	i_{w_1, w_2}\colon\cO_{w_1,w_2}\to X/K', \qquad 
	p_{w_1, w_2}\colon\cO_{w_1, w_2} \to \left(X/K'\right)^{\mathbb{G}_m}\!\!.
	\]
	Since $\Hfb \circ \For^{K'}$ is triangulated, it suffices to check that it does not send nonzero objects to zero. Let $\F'\not\simeq 0$ be a complex in $D^b_{\mathcal L}(K'\backslash X/K')$, and let $\F = \For^{K'}\F'$. Assume that $\Hfb(\F)\simeq0$. Since $\F$ is locally constant along $K'$-orbits, it is in particular locally constant along the fibers of $p_{w_1, w_2}$. Since the fibers of $p_{w_1, w_2}$ are isomorphic to an affine space, we have the implication
	\[
	\Hf_{w_1, w_2}^{\mathrm{res}} (\F) = p_{w_1, w_2, !} \circ i_{w_1, w_2}^*(\F)
	\simeq 0 \quad
	\Longrightarrow \quad
	i_{w_1, w_2}^* \F \simeq 0.
	\]
	Since $\F$ admits a filtration whose successive subquotients are $i_{w_1,w_2,!} \circ i_{w_1,w_2}^*\F$, we have that $\F \simeq 0$ itself. This contradicts $\F'\not\simeq 0$, so we proved that the image of a nonzero object is nonzero.
\end{proof}

 Using the lemma above, it is enough to check that the following transformation is an isomorphism of functors from $D^b_{\mathcal L}(K\backslash X/K')$ 
\begin{equation}\label{eq: the natural transform mon}
	\Hfb \circ \For^{K'}(\alpha') \colon \Hfb \circ \For^{K'} \circ \Av^{K'}_{B,!}\circ \For^{K_0}_B\circ\Limg \to
	\Hfb \circ \For^{K'} \circ \Limg \circ\Av^{K'}_{B,!}\circ \For^K_B.
\end{equation}
The proof repeats the argument in Subsection \ref{subsec:proof_flag} with Proposition~\ref{prop: radon commutes with H monodromic} in place of Proposition~\ref{prop: radon commutes with H}. This establishes Theorem~\ref{sec:character-sheaves-1}.

\subsection{From monodromic sheaves to character sheaves.} \label{subsec:monodrom_to_char}
Here we prove the implication Theorem \ref{sec:character-sheaves-1} $\Rightarrow$ Theorem \ref{sec:character-sheaves}.
The following lemma is immediate.
\begin{lemma}\label{lemma:thick_completion}
	Let $\Phi$, $\Psi \colon\mathcal{C} \to \mathcal{D}$ be two triangulated functors between triangulated categories $\mathcal{C}$ and $\mathcal{D}$, and let $\alpha \colon \Phi \to \Psi$ be a morphism of triangulated functors. Let $C_{ob}$ be a collection of objects in $\mathcal{C}$, and $\mathcal{C}'$ the thick closure of $C_{ob}$.
	
	Assume that $\alpha_c\colon\Phi(c)\to\Psi(c)$ is an isomorphism for any $c\in C_{ob}$. Then the restriction of $\alpha$ to $\mathcal{C}'$ is an isomorphism.
\end{lemma}

Recall that the category of character sheaves $D^b_{\CS_{\mathcal L}}(G)$ is the thick closure of the essential image of $D^b_{\mathcal L}(K\backslash X/K')$ under the functor $\chi$. By Lemma \ref{lemma:thick_completion}, it is enough to show that the natural transformation
\[
	\alpha \circ \chi \colon \Av^{K'}_{B,!}\circ \For^{K_0}_B\circ \Limg \circ \chi \to \Limg \circ\Av^{K'}_{B,!}\circ \For^K_B \circ \chi
\]
restricted to $D^b_{\mathcal{L}}(K \backslash X/K')$ is an isomorphism. To this end, consider the following cube diagram, where parallel arrows in the cube have the same label:
\[
	\definecolor{cof}{RGB}{219,144,71}
	\definecolor{pur}{RGB}{186,146,162}
	\definecolor{greeo}{RGB}{52,111,72}
	\definecolor{greet}{RGB}{52,111,72}
	\definecolor{red-top}{RGB}{186,146,162}
	\definecolor{red-back}{RGB}{140,120,100}
	\definecolor{red-right}{RGB}{160,133,140}
	\definecolor{green-front}{RGB}{91,173,69}
	\definecolor{green-left}{RGB}{52,111,72}
	\definecolor{green-bottom}{RGB}{73,140,70}
	\tikzset{FSTYLE/.style = {draw=blue, % just to demonstrate, where LA is used
			line width=#1, -{Straight Barb[length=3pt]}},
			FSTYLE/.default=1pt
			}
	\tikzset{BSTYLE/.style = {draw=red, % just to demonstrate, where LA is used
			line width=#1, -{Straight Barb[length=3pt]}},
			BSTYLE/.default=1pt
			}
	\begin{tikzcd}
		% [execute at end picture={
		%   \foreach \Nombre in  {A,B,C,D,E,F,G,H,I,J}
		%   {\coordinate (\Nombre) at (\Nombre.center);}
		%   \fill[red-back,opacity=0.2] 
		%   (B) -- (C) -- (G) -- (F) -- cycle;
		%   \fill[green-left,opacity=0.2] 
		%   (A) -- (B) -- (F) -- (E) -- cycle;
		%   \fill[green-bottom,opacity=0.2] 
		%   (E) -- (F) -- (G) -- (H) -- cycle;
		%   \fill[red-top,opacity=0.2] 
		%   (A) -- (B) -- (C) -- (D) -- cycle;
		%   \fill[green-front,opacity=0.2] 
		%   (A) -- (D) -- (H) -- (E) -- cycle;
		%   \fill[red-right,opacity=0.2] 
		%   (H) -- (D) -- (C) -- (G) -- cycle;
		% }]
		&
		|[alias=B]|D^b(K\backslash X/K')
		\ar[FSTYLE, dashed]{dl}[swap]{\Limg}
		\ar{rr}{\chi = \prescript{R}{}{\Av}^{K}_{B,!}\circ \prescript{R}{}{\For}^{K'}_B}
		\ar[BSTYLE]{dd}[near end]{}
		& & |[alias=C]|D^b(K\backslash X/ K)
		\ar{dd}[sloped]{\Av^{K'}_{B,!}\circ \For^K_B}
		\ar{dl}[swap, sloped, near start]{}
		\\
		|[alias=A]|D^b(K_0\backslash X/K')
		\ar[crossing over, FSTYLE, dashed]{rr}[near start]{}
		\ar{dd}[swap, sloped]{\Av^{K'}_{B,!}\circ \For^{K_{0}}_B}
		& |[alias=I]| & |[alias=D]|D^b(K_0 \backslash X/K)
		\\
		&
		|[alias=F]|D^b(K'\backslash X/K')
		\ar[BSTYLE, near start]{rr}{}
		\ar[sloped, near end]{dl}{}
		& |[alias=J]| & |[alias=G]|D^b(K'\backslash X/K)
		\ar[BSTYLE]{dl}{\Limg}
		\\
		|[alias=E]|D^b(K'\backslash X/K')
		\ar{rr}[swap]{\prescript{R}{}{\Av}^{K}_{B,!}\circ \prescript{R}{}{\For}^{K'}_B}
		& & |[alias=H]|D^b(K'\backslash X/K)
		\ar[crossing over, FSTYLE, dashed, from=uu]{}{}
	\end{tikzcd}
\]
Applying Lemma~\ref{lemma:limit-functors3} to the cube, we get that the two natural transformations between the red and the dashed blue paths $D^b(K\backslash X/K') \to D^b(K'\backslash X/K)$ coincide: one given by the composition of natural transformations on the top, right and back faces of the cube, and another along front, left, and bottom faces. Note that the natural transformation corresponding to the right face is exactly the natural transformation $\alpha \circ \chi$. 

The back and front faces are compositions of base change isomorphisms. Since $\chi$ is a composition of smooth pullback and proper pushforward, the top and bottom faces are isomorphisms by Lemma \ref{lemma:limit-functors1}\ref{item:1} and \ref{item:2}. After restricting to the following subcategories, the natural transformation corresponding to the left face is an isomorphism by Theorem~\ref{sec:character-sheaves-1}:
\[
\definecolor{cof}{RGB}{219,144,71}
\definecolor{pur}{RGB}{186,146,162}
\definecolor{greeo}{RGB}{52,111,72}
\definecolor{greet}{RGB}{52,111,72}
\definecolor{red-top}{RGB}{186,146,162}
\definecolor{red-back}{RGB}{140,120,100}
\definecolor{red-right}{RGB}{160,133,140}
\definecolor{green-front}{RGB}{91,173,69}
\definecolor{green-left}{RGB}{52,111,72}
\definecolor{green-bottom}{RGB}{73,140,70}
\tikzset{FSTYLE/.style = {draw=blue, % just to demonstrate, where LA is used
		line width=#1, -{Straight Barb[length=3pt]}},
	FSTYLE/.default=1pt
}
\tikzset{BSTYLE/.style = {draw=red, % just to demonstrate, where LA is used
		line width=#1, -{Straight Barb[length=3pt]}},
	BSTYLE/.default=1pt
}
\begin{tikzcd}
	% [execute at end picture={
		%   \foreach \Nombre in  {A,B,C,D,E,F,G,H,I,J}
		%   {\coordinate (\Nombre) at (\Nombre.center);}
		%   \fill[red-back,opacity=0.2] 
		%   (B) -- (C) -- (G) -- (F) -- cycle;
		%   \fill[green-left,opacity=0.2] 
		%   (A) -- (B) -- (F) -- (E) -- cycle;
		%   \fill[green-bottom,opacity=0.2] 
		%   (E) -- (F) -- (G) -- (H) -- cycle;
		%   \fill[red-top,opacity=0.2] 
		%   (A) -- (B) -- (C) -- (D) -- cycle;
		%   \fill[green-front,opacity=0.2] 
		%   (A) -- (D) -- (H) -- (E) -- cycle;
		%   \fill[red-right,opacity=0.2] 
		%   (H) -- (D) -- (C) -- (G) -- cycle;
		% }]
	&
	|[alias=B]|D^b_{\mathcal{L}}(K\backslash X/K')
	\ar[FSTYLE, dashed]{dl}[swap]{\Limg}
	\ar{rr}{\chi = \prescript{R}{}{\Av}^{K}_{B,!}\circ \prescript{R}{}{\For}^{K'}_B}
	\ar[BSTYLE]{dd}[near end]{}
	& & |[alias=C]|D^b_{\CS_{\mathcal{L}}}(G)
	\ar{dd}[sloped]{\Av^{K'}_{B,!}\circ \For^K_B}
	\ar{dl}[swap, sloped, near start]{}
	\\
	|[alias=A]|D^b(K_0\backslash X/K')
	\ar[crossing over, FSTYLE, dashed]{rr}[near start]{}
	\ar{dd}[swap, sloped]{\Av^{K'}_{B,!}\circ \For^{K_0}_B}
	& |[alias=I]| & |[alias=D]|D^b(K_{0} \backslash X/K)
	\\
	&
	|[alias=F]|D^b(K'\backslash X/K')
	\ar[BSTYLE, near start]{rr}{}
	\ar[sloped, near end]{dl}{}
	& |[alias=J]| & |[alias=G]|D^b(K'\backslash X/K)
	\ar[BSTYLE]{dl}{\Limg}
	\\
	|[alias=E]|D^b(K'\backslash X/K')
	\ar{rr}[swap]{\prescript{R}{}{\Av}^{K}_{B,!}\circ \prescript{R}{}{\For}^{K'}_B}
	& & |[alias=H]|D^b(K'\backslash X/K)
	\ar[crossing over, FSTYLE, dashed, from=uu]{}{}
\end{tikzcd}
\]
For the latter cube, the natural transformations obtained by composing along the top, right, and back faces and along the front, left, and bottom faces coincide. We have established that five out of the six transformations are isomorphisms, so the remaining one is an isomorphism as well. That is, $\alpha \circ \chi$ is an isomorphism, which implies Theorem~\ref{sec:character-sheaves}.

\section{Generalization to a pair of arbitrary parabolic subgroups}\label{sec:gener-pair-arbitr}
\subsection{Limit functors.}
For a parabolic subgroup $P$ and its opposite $P^-$ with common Levi $L$, denote
\[
K_P = P \times_L P^-, \qquad K'_P = P \times_L P.
\]
Choose a cocharacter $\lambda_P\colon\mathbb{G}_m \to T$ such that $P$ and $P^-$ are respectively the attracting and repelling loci for the adjoint $\mathbb{G}_m$-action on $G$ defined by $\lambda_P$. The image of $\lambda_P$ is in the center of $L = P \cap P^-$. Define a cocharacter $\gamma_P$ of $G \times G$ by $\gamma_P = (\lambda_P^{-1}, \lambda_P)$. Consider another parabolic subgroup $Q\supset P$. We have a group scheme 
\[
	\widetilde{K}^Q_P \subset G \times G \times \mathbb{P}^1
\] 
over $\mathbb{P}^1$ whose fiber over $t \in \mathbb{G}_m \subset \mathbb{P}^1$ is $\gamma_P(t)K_Q\gamma_P(t)^{-1}$ and whose fiber over $0$ is $K_P$.
As before, this group scheme gives rise to the limit functor
\[
\Lim_{\gamma_P}\colon D^b(K_Q \backslash X) \to D^b(K_P \backslash X)
\]
for a scheme $X$ with a left action of $G \times G$.

\begin{remark}
	Subsection~\ref{paragraph: limit functor} corresponds to the specialization of this for 
	$Q=G$ and $P=B$. In this case, $K_P$, $K'_P$, and $K_Q$ become $K_B = K_0$, $K'_B = K'$, and $K_G = K$, respectively.
\end{remark}

\subsection{Hyperbolic localization.}
For a fixed parabolic subgroup $P$ containing $B$, let $W_P$ be the corresponding parabolic subgroup of $W$. Consider the following action of $\mathbb{G}_m$ on $X/K'_B$, for $w_1, w_2 \in W_P\backslash W$:
\begin{gather*}
	\aw^P\colon \mathbb{G}_m \times X/K'_B \rightarrow X/K'_B\\ 
	\aw^P\left(t, (g_1, g_2) \right) = \Big( \lambda_P(t) \cdot g_1 \cdot w_1^{-1} \left( \lambda_P (t^{-1}) \right),
	\lambda_P(t) \cdot g_2 \cdot w_2^{-1} \left(\lambda_P \left(t^{-1} \right) \right)\Big).
\end{gather*}
This gives us the hyperbolic localization functors
\[
	\Hf_{w_1, w_2}^P\colon D^b_{mon}(X/K'_B) \rightarrow  D^b((X/K'_B)^{\mathbb{G}_m}).
\]
The following lemma follows directly from the parabolic Bruhat decomposition.
\begin{lemma} 
	For the $\mathbb{G}_m$-action on $X/K'_B$ given by $\aw^P$, the following hold.
	\begin{enumerate}[label=(\alph*),ref=(\alph*)]
		\item The fixed locus 
		\[
			(X/K'_B)^{\mathbb{G}_m} = \bigsqcup_{\tilde w_1,\tilde w_2} (L \tilde{w}_1 B \times L \tilde{w}_2 B)/K'_B,
		\]
		where the disjoint union is taken over the pairs $\tilde{w}_1,\tilde{w}_2 \in W_P \backslash W$ satisfying
		\begin{equation}%\label{eq:involved_condition}
			w_1^{-1}(\lambda_P) - w_2^{-1}(\lambda_P) = \tilde{w}_1^{-1}(\lambda_P) - \tilde{w}_2^{-1}(\lambda_P).
		\end{equation}
		\item The attracting locus for $(L \tilde{w}_1 B \times L \tilde{w}_2 B)/K'_B$ is $(P \tilde{w}_1 B \times P\tilde{w}_2 B)/K'_B$.
	\end{enumerate}
\end{lemma}

Consider the restriction of $\Hf_{w_1, w_2}^P$
\[
	\Hf^{P, \mathrm{res}}_{w_1, w_2}\colon D^b_{\mathcal L}(X/K'_B) \to D^b \big( Lw_1B\times Lw_2B/K'_B\big)
\]
where the subscript $\mathcal L$ stands for $\mathcal L$-monodromic for the right $T$-action. We denote
\[
	\Hfb^P = \bigoplus_{w_1, w_2 \in W_P\backslash W} \Hf_{w_1, w_2}^{P, \mathrm{res}}.
\]

We have the following parabolic generalization of Lemma \ref{sec:proof-monodr-sheav}; the proof is completely analogous.
\begin{lemma}\label{sec:hyperb-local}
		The functor 
		\[
		\Hf_{\bullet}^P \circ \For^{K'_P}\colon D^b_{\mathcal L}(K'_P\backslash X/K'_B) \rightarrow \bigoplus_{w_1, w_2 \in W_P \backslash W} D^b \big( Lw_1B \times Lw_2B/K'_B \big)
		\]
		is conservative.
\end{lemma}

\subsection{Parabolic character sheaves.}

In this subsection, we recall the definition of parabolic character sheaves used in \cite{lusztigParabolicCharacterSheaves2003} and compare it with the characterization needed here. Alternatively, one may take Proposition~\ref{prop:non-prime_parabolic} as the definition of the block $D^b_{\CS_{W\!\mathcal L}}(K_P\backslash X/K_G)$, and define the category of parabolic character sheaves $D^b_{\CS}(K_P\backslash X/K_G)$ to be generated by these blocks as $\mathcal L$ ranges over rational rank-one local systems. The second mentioned definition is parallel to the one used in \cite{lusztigParabolicCharacterSheaves2003a}.

The following additional assumption is needed for the comparison of definitions. Recall the sheaf-theoretic settings of the Introduction. In the case of sheaves with coefficients in a field $\mathbb{k}$ on stacks over $\mathbb{C}$, we will assume that $\operatorname{char} \mathbb k$ does not divide $|W|$.

As above, we consider the scheme $X = G \times G$ with the actions of $G\times G$ by left and right multiplication.

We define the following parabolic analogue of the Harish-Chandra functor:
\[
	\hc^Q_P \colon D^b(K_{Q}'\backslash X/K_G) \to D^b(K_{P}'\backslash X/K_G), \quad \hc^Q_P = {\Av}^{K'_P}_{K'_Q \cap K'_P,!}\circ{\For}^{K'_Q}_{K'_Q\cap K'_P}.
\]
Also consider the corresponding character functor
\[
	\chi^Q_P \colon D^b(K'_P\backslash X/K_G) \to D^b(K_{Q}'\backslash X/K_G), \quad \chi^Q_P = {\Av}^{K'_Q}_{K'_P \cap K'_Q,!}\circ{\For}^{K'_P}_{K'_P\cap K'_Q}.
\]
Note that $\chi^Q_P$ is right adjoint to $\hc^Q_P$ up to a homological shift, since $\chi^Q_P$ is a composition of smooth pullback and proper pushforward, the latter corresponding to a map whose fibers are isomorphic to $K_Q'/(K'_P\cap K'_Q) \simeq Q/P$.

The functors $\hc_P^Q$ and $\chi_P^Q$ satisfy the following transitivity property.
\begin{proposition}[{\cite[Lemma~6.2]{lusztigParabolicCharacterSheaves2003}}]
\label{sec:parab-char-sheav}
	Let $R \subset P \subset Q$ be a triple of parabolic subgroups. We have the following natural isomorphisms of functors:
	\[
		\hc^P_R \circ \hc^Q_P \simeq \hc^Q_R, \qquad \chi^Q_P \circ \chi^P_R \simeq \chi^Q_R.
	\]
\end{proposition}

We now recall the notion of \emph{parabolic character sheaves} from \cite{lusztigParabolicCharacterSheaves2003}.
The torus $T$ acts on $K'_B\backslash X/K_G$ on the left. Let $\mathcal{L}$ be a rational rank-one local system on $T$. The derived category of parabolic character sheaves with character $\mathcal{L}$ is defined as the thick subcategory of $D^b(K'_P\backslash X/K_G)$ generated by the essential image of $D^b_{\mathcal L}(K'_B\backslash X/K_G)$ under the functor $\chi_B^P$; it will be denoted by $D^b_{\CS_{\mathcal L}}(K'_P\backslash X/K_G)$. We write $D^b_{\CS}(K'_P\backslash X/K_G)$ for the full triangulated subcategory of $D^b(K'_P\backslash X/K_G)$ generated by $D^b_{\CS_{\mathcal{L}}}(K'_P\backslash X/K_G)$ for all rational $\mathcal{L}$.

\begin{lemma}[{\cite[Proposition 6.7]{lusztigParabolicCharacterSheaves2003}}]\label{sec:parab-char-sheav-1}
	The functor $\hc^Q_P$ sends $D^b_{\CS}(K_Q'\backslash X/K_G)$ to $D^b_{\CS}(K'_P\backslash X/K_G).$ The functor $\chi^Q_P$ sends $D^b_{\CS}(K'_P\backslash X/K_G)$ to $D^b_{\CS}(K_Q'\backslash X/K_G).$
\end{lemma}

\begin{remark}
	Note that the pair of functors $(\hc^G_B, \chi^G_B)$ differs from the previously defined pair of functors $(\hc, \chi)$ by interchanging the roles of the left and right actions. This does not affect the resulting subcategory $D^b_{\CS_{\mathcal{L}}}(K'_G\backslash X/K_G) = D^b_{\CS_{\mathcal{L}}}(G)$, as we will see from the following discussion. In the setting of sheaves with coefficients in $\mathbb k$, this conclusion uses our assumption on $\operatorname{char} \mathbb k$.
\end{remark}

We use the notation $D^b_{W\!\mathcal L}$ for the category generated by $w\mathcal L$-monodromic sheaves as $w$ ranges over $W$, and call its objects $W\mathcal L$-monodromic sheaves.

Parabolic character sheaves $D^b_{\CS_{\mathcal{L}}}(K'_P\backslash X/K_G)$ coincide with $D^b_{\CS_{\mathcal{L'}}}(K'_P\backslash X/K_G)$ if and only if $\mathcal L'$ is in the $W_P$-orbit of $\mathcal L$. We therefore write $D^b_{\CS_{W\!\mathcal L}}(K'_P\backslash X/K_G)$ for the subcategory generated by $D^b_{\CS_{w\mathcal L}}(K'_P\backslash X/K_G)$ for $w\in W$.

\begin{proposition}\label{prop:parabolic_char_right_left}
	The category $D^b_{\CS_{W\!\mathcal L}}(K'_P\backslash X/K_G)$ is the idempotent completion of the essential image of $D^b_{W\!\mathcal L}(K'_P\backslash X/K'_B)$ under $\prescript{R}{}{\Av}^{K_G}_{K'_B\cap K_G,!} \circ \prescript{R}{}{\For}^{K'_B}_{K'_B\cap K_G}$, where $D^b_{W\mathcal{L}}(K'_P \backslash X/K'_B)$ is the $W\mathcal{L}$-monodromic category for the right $T$-action.
\end{proposition}

We give the proof of the proposition in Appendix~\ref{appendix:proof_left_right}.

\bigskip
We now pass from $K'_P\backslash X/K_G$ to $K_P\backslash X/K_G$.
To define parabolic character sheaves in the latter setting, we need the following proposition.

\begin{proposition}\label{prop:parab-char-sheav-2}
	The following version of the parabolic Radon transform functor 
	\begin{equation} \label{eq:parabolic_radon}
		{\Av}^{K_P}_{K_P \cap K'_P,!}\circ{\For}^{K'_P}_{K_P\cap K'_P}\colon D^b_{\CS}(K'_P\backslash X/K_G) \to D^b(K_P\backslash X/K_G)
	\end{equation}
	is fully faithful.
\end{proposition}
We give the proof of the proposition in Appendix~\ref{appendix:fully_faith_R}.

\bigskip

For each rational $\mathcal L$, we define the category $D^b_{\CS_{\mathcal L}}(K_P\backslash X/K_G)$ of parabolic character sheaves on $K_P\backslash X/K_G$ as the essential image of $D^b_{\CS_{\mathcal L}}(K'_P\backslash X/K_G)$ under the parabolic Radon transform~\eqref{eq:parabolic_radon}. The preceding discussion shows that this category is thick. Let $D^b_{\CS}(K_P\backslash X/K_G)$ be the full triangulated subcategory of $D^b(K_P\backslash X/K_G)$ generated by $D^b_{\CS_{\mathcal{L}}}(K_P\backslash X/K_G)$ for all rational $\mathcal{L}$.

Below we use that $D^b_{\CS}(K_P\backslash X/K_G)$ is generated by the blocks $D^b_{\CS_{W\!\mathcal L}}(K_P\backslash X/K_G)$. It follows from the above that the latter subcategories admit the following characterization.

\begin{proposition}\label{prop:non-prime_parabolic}
	The category $D^b_{\CS_{W\!\mathcal{L}}}(K_P\backslash X/K_G)$ can be described as the thick subcategory of $D^b(K_P\backslash X/K_G)$ generated by the essential image of $D^b_{W\mathcal{L}}(K_P \backslash X/K'_B)$ under the functor $\prescript{R}{}{\Av}^{K_G}_{K'_B \cap K_G,!}\circ\prescript{R}{}{\For}^{K'_B}_{K'_B\cap K_G}$.
\end{proposition} 

\subsection{The main result in the parabolic setting.} \label{sec:main_res_parabol}
Let us now fix $P$, $Q$, and $\gamma_P$ for the rest of this subsection. To keep our notation light, we will write $\gamma$ for $\gamma_P$.

Consider the diagram:
\begin{equation*}  %\label{eq:diagram_main_parabolic}
	\begin{tikzcd}
		D^b(K_Q\backslash X/K_G) \arrow[r, "\Limg"]\arrow[d, "\Av^{K'_P}_{P,!}\circ \For^{K_Q}_P"'] & D^b(K_P\backslash X/K_G) \arrow[d,"\Av^{K'_P}_{P,!}\circ \For^{K_P}_P"] \arrow[dl, Rightarrow, "\alpha_P^Q"]\\
		D^b(K'_P\backslash X/K_G) \arrow[r, "\Limg"'] & D^b(K'_P\backslash X/K_G)
	\end{tikzcd}
\end{equation*}

Here $\alpha_P^Q \colon \Av^{K'_P}_{P,!}\circ \For^{K_P}_P\circ\Limg \to \Limg\circ \Av^{K'_P}_{P,!}\circ \For^{K_Q}_P$ is the natural transformation given by Lemma \ref{lemma:limit-functors1}.

\begin{theorem}\label{sec:main-result-parab}
	The restriction of $\alpha_P^Q$ to the subcategory of parabolic character sheaves $D^b_{\CS}(K_Q\backslash X/K_G)$ is an isomorphism.
\end{theorem}

Recall that in Subsection \ref{subsec:monodrom_to_char} we derived Theorem \ref{sec:character-sheaves} from Theorem \ref{sec:character-sheaves-1}. Using a similar argument together with the description of $D^b_{\CS_{W\!\mathcal{L}}}(K_Q\backslash X/K_G)$ given by Proposition~\ref{prop:non-prime_parabolic}, one can reduce Theorem \ref{sec:main-result-parab} to the following statement for monodromic sheaves. Consider the diagram:
\begin{equation*} %\label{eq:parabolic-monodromic}
	\begin{tikzcd}
		D^b(K_Q\backslash X/K'_B) \arrow[r, "\Limg"]\arrow[d, "\Av^{K'_P}_{P,!}\circ \For^{K_Q}_P"'] & D^b(K_P\backslash X/K'_B) \arrow[d,"\Av^{K'_P}_{P,!}\circ \For^{K_P}_P"] \arrow[dl, Rightarrow, "\tilde{\alpha}^Q_P"]\\
		D^b(K'_P\backslash X/K'_B) \arrow[r, "\Limg"'] & D^b(K'_P\backslash X/K'_B)
	\end{tikzcd}
\end{equation*}
Here $\tilde{\alpha}^Q_P \colon \Av^{K'_P}_{P,!}\circ \For^{K_P}_P\circ\Limg \to \Limg\circ\Av^{K'_P}_{P,!}\circ \For^{K_Q}_P$ is again the natural transformation given by Lemma \ref{lemma:limit-functors1}. 

\begin{theorem}\label{sec:parabolic-monodromic-sheaves}
	For any rational rank-one local system $\mathcal{L}$ on $T$, the restriction of $\tilde{\alpha}^Q_P$ to the monodromic category $D^b_{\mathcal{L}}(K_Q\backslash X/K'_B)$ is an isomorphism.
\end{theorem}

The proof is identical to that of Theorem~\ref{sec:character-sheaves-1}, with Lemma~\ref{sec:hyperb-local} in place of Lemma~\ref{sec:proof-monodr-sheav}.

Consider the following versions of the parabolic Radon transform
\begin{align*}
	R_{!}^P \colon& D^b_{\CS}(K_{P}\backslash X/K_G) \to D^b_{\CS}(K'_{P}\backslash X/K_G), \quad R_{!}^P = \Av^{K'_P}_{P,!}\circ \For^{K_P}_P.	\\
	R_{*}^P \colon& D^b_{\CS}(K'_{P}\backslash X/K_G) \to D^b_{\CS}(K_{P}\backslash X/K_G), \quad R_{*}^P = \Av^{K_P}_{P,*}\circ \For^{K'_P}_P.	
\end{align*}
The Verdier-dual version of Proposition \ref{prop:parab-char-sheav-2} implies that $R_*^P$ is invertible. Since $R_!^P$ is the left adjoint of $R_*^P$, the functors are mutually inverse.

\begin{corollary}
	There is an isomorphism of functors on $D^b_{\CS}(K_Q\backslash X/K_G)$
	\[
		\Limg \simeq R^P_* \circ \hc^Q_P \circ R^Q_!.
	\]
\end{corollary}

The proof is similar to that of Corollary \ref{RHC}; the only
additional ingredient is the following base-change isomorphism:
\[
	\Av^{K'_P}_{P,!}\circ \For^{K_Q}_P \simeq \hc^Q_P \circ R^Q_!.
\]

\iffalse
\begin{proof}
	By Lemma \ref{lemma:limit-functors2}, the lower horizontal arrow in the diagram \eqref{eq:diagram_main_parabolic} is isomorphic to the identity functor when restricted to the monodromic subcategory. We obtain an isomorphism
	\[
		\Av^{K'_P}_{P,!}\circ \For^{K_Q}_P \simeq \Av^{K'_P}_{P,!}\circ \For^{K_P}_P	\circ \Limg
	\]
	\[
		\Limg \simeq \Av^{K_P}_{P,*}\circ \For^{K'_P}_P \circ
		\Av^{K'_P}_{P,!}\circ \For^{K_Q}_P
	\]
	Proposition \ref{sec:hyperb-local-parab} implies that $\Av^{K_P}_{P,*}\circ \For^{K'_P}_P$ is the inverse of $\Av^{K'_P}_{P,!}\circ \For^{K_P}_P$.
\end{proof}
\fi

\appendix

\section{On parabolic character sheaves}

In this appendix, we prove Propositions~\ref{prop:parabolic_char_right_left} and~\ref{prop:parab-char-sheav-2}. These results are needed to compare the definition of parabolic character sheaves used in \cite{lusztigParabolicCharacterSheaves2003} with the characterization used in Subsection~\ref{sec:main_res_parabol}. In the setting of sheaves with coefficients in $\mathbb k$, the results below use our
assumption that $\operatorname{char}\mathbb k$ does not divide $|W|$.

\subsection{Proof of Proposition \ref{prop:parabolic_char_right_left}} \label{appendix:proof_left_right}

\begin{lemma}\label{sec:parab-radon-transf}
	Let $X$ be a scheme with a left action of $G \times G$, and a right action of a group $H$ commuting with the action of $G \times G$.
	Let $P \subset Q$ be a pair of parabolic subgroups. Then the composition 
	\begin{equation} \label{eq:compos_chi_hc}
		{}{\Av}^{K'_Q}_{K'_P \cap K'_Q,!}\circ{}{\For}^{K'_P}_{K'_P\cap K'_Q} \circ {}{\Av}^{K'_P}_{K'_P \cap K'_Q,!}\circ{}{\For}^{K'_Q}_{K'_P\cap K'_Q}
	\end{equation}
	contains a homological shift of the identity functor on $D^b(K_Q'\backslash X/H)$ as a direct summand. The analogous statement holds for $X$ with the right action of $G\times G$.
\end{lemma}

\begin{proof}
	A similar proof of a less general version can be found, for example, in \cite[Proposition 3.3.1]{tolmachovLinearStructureFinite2024}; see also \cite[Lemma 6.6]{lusztigParabolicCharacterSheaves2003}. For the convenience of the reader, we reproduce the proof using our notation. Note that the category $D^b(K_Q'\backslash X/H)$ admits a left action of the monoidal category $D^b(K_Q' \backslash (G \times G)/ K_Q')$, with both the monoidal and the module operations being given by the usual $!$-convolution. Consider the map \[\mu\colon (K'_P\cap K_Q')\backslash K'_P / (K'_P\cap K_Q') \to K_Q' \backslash (G \times G)/K_Q'\] induced by the embeddings of the groups. The functor \eqref{eq:compos_chi_hc} can be identified with convolution with the sheaf $\mu_!\cc$. It follows from Springer theory for the Levi of $Q$ \cite{borhoRepresentationsGroupesWeyl1981} that $\mu_!\cc$ contains a homological shift of $\iota_*\cc_{K_Q'\backslash K_Q'/K_Q'}$ as a direct summand, where $\iota\colon K_Q'\backslash K_Q'/K_Q' \to K_Q'\backslash (G\times G)/K_Q'$ is the embedding. This concludes the proof.
\end{proof}
\begin{remark}
	In the setting of sheaves with coefficients in a field $\mathbb k$, we used above that the trivial representation of the Weyl group $W_Q$ is a direct summand of the permutation representation $\mathbb{k}[W_Q/W_P]$, as follows from our assumption on $\operatorname{char} \mathbb k$. By Springer theory, this representation-theoretic decomposition induces the corresponding decomposition of the parabolic Springer sheaf $\mu_!\cc$.
\end{remark}

\begin{corollary}\label{corol:thick_vs_idemp}
	The category $D^b_{\CS_{\mathcal{L}}}(K'_P\backslash X/K_G)$ is the idempotent completion of the essential image of $D^b_{\mathcal L}(K'_B\backslash X/K_G)$ under the functor $\chi_B^P$.
\end{corollary}
\begin{proof}
	Recall that we defined $D^b_{\CS_{\mathcal{L}}}(K'_P\backslash X/K_G)$ as the thick closure of the essential image of $\chi_B^P$, so we need to check that the idempotent completion in question is a full triangulated subcategory. This follows straightforwardly from Lemma~\ref{sec:parab-radon-transf}.
\end{proof}

\begin{proof}[Proof of Proposition \ref{prop:parabolic_char_right_left}]
	%We prove the proposition by showing that 
	It is sufficient to show that
	both subcategories coincide with the idempotent completion of the essential image of $D^b_{W\mathcal{L},W\mathcal{L}}(K'_B \backslash X/K'_B)$ under the functor
	\[
	{\Av}^{K'_P}_{K'_P \cap K'_B,!}\circ{\For}^{K'_B}_{K'_B\cap K'_P} \circ \prescript{R}{}{\Av}^{K_G}_{K'_B \cap K_G,!}\circ\prescript{R}{}{\For}^{K'_B}_{K'_B\cap K_G},
	\] 
	where $D^b_{W\mathcal L,W\mathcal L}(K'_B\backslash X/K'_B)$ denotes the category of sheaves that are $W\mathcal L$-monodromic for both the left and right $T$-actions.
	
	%First note that the essential image of $D^b_{{\mathcal{L}}}(K'_B\backslash X/K)$ in $D^b_{\CS_{\mathcal{L}}}(K'_P\backslash X/K)$ coincides with the essential image of $D^b_{W_P\mathcal{L}}(K'_B\backslash X/K)$, where the subscript $W_P\mathcal{L}$ again stands for all the monodromies in the $W_P$-orbit of $\mathcal{L}$.
	\begin{lemma}\label{lemma:monodromic_W_P_W}
		The image of $D^b_{W\mathcal{L}}(K'_B \backslash X/K_G)$ under $\prescript{R}{}{\Av}^{K'_B}_{K'_B\cap K_G, !} \circ \prescript{R}{}{\For}^{K_G}_{K'_B\cap K_G}$ is contained in $D^b_{W\mathcal{L},W\mathcal{L}}(K'_B \backslash X/K'_B)$.
	\end{lemma}
	\begin{proof}
		For $w_1,w_2\in W$, let $i_{w_1,w_2}$ be the inclusion of the Bruhat cell indexed by $(w_1,w_2)$. Any object $\mathcal{F} \in D^b(K'_B \backslash X/K'_B)$ admits a filtration whose successive subquotients are given by $i_{w_1,w_2,*} \circ i_{w_1,w_2}^!\mathcal{F}$. 
		
		The functor in the lemma statement preserves $W\mathcal L$-monodromic
		sheaves. Therefore, if $\mathcal F$ belongs to its essential image, then the
		subquotients are $(W\mathcal L,W\mathcal L)$-monodromic. Hence $\mathcal F$
		is also $(W\mathcal L,W\mathcal L)$-monodromic.
	\end{proof}
	\begin{lemma}
		The category $D^b_{W\mathcal{L}}(K'_B \backslash X/K_G)$ is the idempotent completion of the essential image of $D^b_{W\mathcal{L},W\mathcal{L}}(K'_B \backslash X/K'_B)$ under the functor $\prescript{R}{}{\Av}^{K_G}_{K'_B \cap K_G,!}\circ\prescript{R}{}{\For}^{K'_B}_{K'_B\cap K_G}$.
	\end{lemma}
	\begin{proof}
		Lemma \ref{lemma:monodromic_W_P_W} and Lemma~\ref{sec:parab-radon-transf} imply that $D^b_{W\mathcal{L}}(K'_B \backslash X/K_G)$ is contained in the idempotent completion of the essential image. On the other hand, the functor evidently sends left-monodromic complexes to left-monodromic complexes with the same monodromy.
	\end{proof}

	Similarly, the category $D^b_{W\mathcal{L}}(K'_P \backslash X/K'_B)$ is the idempotent completion of the essential image of $D^b_{W\mathcal{L},W\mathcal{L}}(K'_B \backslash X/K'_B)$ under the functor ${\Av}^{K'_P}_{K'_P \cap K'_B,!} \circ {\For}^{K'_B}_{K'_B\cap K'_P}$. Then the proposition follows from the base-change isomorphism.
\end{proof}

\subsection{Proof of Proposition \ref{prop:parab-char-sheav-2}.} \label{appendix:fully_faith_R}
To prove the proposition, we need the following result.

\begin{proposition}\label{sec:hyperb-local-parab}
	The following version of the parabolic Radon transform functor 
	\[
	{\Av}^{K_P}_{K_P \cap K'_P,!}\circ{\For}^{K'_P}_{K_P\cap K'_P}\colon D^b_{\mathcal{L}}(K'_P\backslash X/K'_B) \to D^b_{\mathcal{L}}(K_P\backslash X/K'_B)
	\]
	is an equivalence of categories.
\end{proposition}
\begin{proof}
	The proof is the same as that of \cite[Theorem 5.2]{chenFormulaGeometricJacquet2017}, using Lemma \ref{sec:hyperb-local}.
\end{proof}

\begin{proof}[Proof of Proposition \ref{prop:parab-char-sheav-2}.]
	The essential image of $D^b_{\CS_{\mathcal{L}}}(K'_P\backslash X/K_G)$ under the functor $\prescript{R}{}{\Av}^{K'_B}_{K'_B\cap K_G, !}\circ\prescript{R}{}{\For}^{K_G}_{K'_B\cap K_G}$ lies in $D^b_{\mathcal{L}}(K'_P\backslash X/K'_B)$. The base-change isomorphism implies that the following diagram commutes:
	\begin{equation*}
		\begin{tikzcd}
			D^b_{\CS_{\mathcal{L}}}(K'_P\backslash X/K_G) \arrow[r]\arrow[d] & D^b(K_P\backslash X/K_G)\arrow[d]\\
			D^b_{\mathcal{L}}(K'_P\backslash X/K'_B) \arrow[r]& D^b(K_P\backslash X/K'_B)
		\end{tikzcd}
	\end{equation*}
	where both horizontal arrows are given by ${\Av}^{K_P}_{K_P \cap K'_P,!}\circ{\For}^{K'_P}_{K_P\cap K'_P}$ and both vertical arrows are given by $\prescript{R}{}{\Av}^{K'_B}_{K'_B\cap K_G, !}\circ\prescript{R}{}{\For}^{K_G}_{K'_B\cap K_G}$.

	Lemma~\ref{sec:parab-radon-transf} implies that the vertical arrows are faithful. Together with Proposition~\ref{sec:hyperb-local-parab}, this shows that the functor~\eqref{eq:parabolic_radon} is faithful. Since $D^b_{\CS_{\mathcal{L}}}(K'_P\backslash X/K_G)$ is the idempotent completion of the essential image of $\chi_B^P$, Lemma~\ref{sec:parab-radon-transf} implies that the functor is also full.
\end{proof}

\bibliography{bibliography}
\bibliographystyle{alpha}

\Addresses

\end{document}